\documentclass[12 pt]{article}
\usepackage{amsfonts}
\title{Domination Cover Pebbling:  Graph Families}
\author{James Gardner\\
Department of Mathematics and Computer Science\\
Emory University\and
Anant P.~Godbole\\
Department of Mathematics\\
East Tennessee State University\and
Alberto Mokak Teguia\\
Department of Mathematics\\
Duke University\and
Annalies Z.~Vuong\\
Department of Mathematics\\
Dartmouth College\and
Nathaniel Watson\\
Department of Mathematics\\
Washington University at St.~Louis\and
Carl R.~Yerger\\
Department of Mathematics\\
Georgia Institute of Technology}
\begin{document}
\def\qed{\vbox{\hrule\hbox{\vrule\kern3pt\vbox{\kern6pt}\kern3pt\vrule}\hrule}}
\def\ep{\varepsilon}
\def\lr{\left(}
\def\lf{\lfloor}
\def\rf{\rfloor}
\def\lc{\left\{}
\def\rc{\right\}}
\def\rr{\right)}
\def\etapn{{2^{n+1} \left({\frac{1-8^{-(k_n+1)}}{7}}\right)+\floor{\frac{\alpha_n}{2}}}}
\def\letapn{{2^{n+1} \left({(1-8^{-(k_n+1)})}/7\right)+\floor{\frac{\alpha_n}{2}}}}
\def\p{\mathbb P}
\def\v{\mathbb V}
\def\e{\mathbb E}
\def\l{\mathbb L}
\def\lg{{\rm lg}}
\newtheorem{thm}{Theorem}
\newtheorem{lem}[thm]{Lemma}
\newtheorem{prop}[thm]{Proposition}
\newtheorem{defn}[thm]{Definition}
\newcommand{\ignore}[1]{}

\providecommand{\abs}[1]{\vert#1\vert}
\providecommand{\norm}[1]{\Vert#1\Vert}
\providecommand{\flooralpha}{\left\lfloor\frac{\alpha_n}{2}\right\rfloor}
\providecommand{\floor}[1]{\left\lfloor#1\right\rfloor}
\maketitle
\begin{abstract}
Given a configuration of pebbles on the vertices of a connected graph $G$, a {\it pebbling move} is defined as the removal of two pebbles from some vertex, and the placement  of one of these on an adjacent vertex.  We introduce the notion of domination cover pebbling, obtained by combining graph cover pebbling (\cite{zsu}) with the theory of domination in graphs (\cite{haynes}). The domination cover pebbling
number, $\psi(G)$, of a graph $G$ is the minimum number of pebbles that must be placed on $V(G)$ such that after a sequence of pebbling moves, the set of
vertices with pebbles forms a dominating set of $G$, regardless of the initial configuration of pebbles.  We discuss basic results and determine $\psi(G)$ for
paths, cycles and complete binary trees.
\end{abstract}

\section{Introduction}
One recent development in graph theory, suggested by Lagarias and Saks and
called pebbling, has been the subject of much research.  It was
first introduced into the literature by Chung~\cite{chung}, and has
been developed by many others including Hurlbert, who published a survey of
pebbling results in~\cite{hurlbert}. There have been many developments since Hurlbert's survey appeared; some of these are described in this paper.

Given a graph $G$, distribute $k$
pebbles (indistinguishable markers) on its vertices in some configuration $C$.  Specifically,
a configuration on a graph $G$ is a function from $V(G)$ to $\mathbb{N}
\cup \{0\} $ representing an arrangement of pebbles on $G$.  For our
purposes, we will always assume that $G$ is connected.  A pebbling move is
defined as the removal of two pebbles from some vertex and the placement of one of these pebbles pebbles on an
adjacent vertex.  Define the pebbling number, $\pi(G)$ to be the minimum number of
pebbles such no matter what their initial configuration, it is possible to move to any root vertex $v$ in $G$ after a sequence of pebbling moves.  Implicit in this definition is the fact that if after
moving to vertex $v$ one desires to move to another root
vertex, the pebbles reset to their original initial configuration.

In this paper, we will combine two ideas, {\it cover} pebbling (\cite{zsu}) and domination (\cite{haynes}), to
introduce a new graph invariant called the domination cover pebbling
(DCP) number of a graph, denoted by $\psi(G)$.  Recall that a set of vertices $D$
in $G$ is a dominating set if every vertex in $G$ is either in $D$ or
adjacent to some element in $D$. The cover pebbling number, $\lambda(G)$ is
defined as the minimum number of pebbles required such that given any
initial configuration of at least $\lambda(G)$ pebbles, it is possible to
make a series of pebbling moves to place at least one pebble on {\it every}
vertex of $G$.  The domination cover pebbling number of a graph $G$, proposed by A.~Teguia,  is the
 minimum number of pebbles required so that any initial
configuration of pebbles can be transformed by a sequence of pebbling moves
so that the set of vertices that contain pebbles form a dominating set $S$
of $G$.  The motivation for our definition comes from a hypothetical situation in which one wishes to transport monitors along the edges of a network that could ultimately ``watch" each vertex -- but half the devices are lost during each move.  The pebbles may be placed on any of the vertices of $G$, and $S$ depends, in general, on the initial configuration -- most importantly, however, $S$ need not equal a minimum dominating set.
Graphs can be easily constructed to illustrate these facts.  Consider the
configurations of pebbles on $P_4$, the path on four vertices, as shown in Figure 1:
\begin{figure}[htb] \unitlength 1mm
\begin{center}
\begin{picture}(100,30)
\put(5,5){\circle*{3}} \put(15,5){\circle*{3}} \put(25,5){\circle*{3}}
\put(35,5){\circle*{3}}\put(55,5){\circle*{3}}
\put(65,5){\circle*{3}}\put(75,5){\circle*{3}} \put(85,5){\circle*{3}}
\put(5,5){\line(1,0){30}} 
\put(55,5){\line(1,0){30}} \put(4,10){5} \put(84,10){5}
\end{picture}
\end{center}
\caption{An example where two different initial configurations produce two
different domination cover solutions.} \label{ex2}
\end{figure}
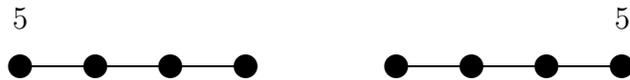
For the graph on the left, we make pebbling moves so that the first and
third vertices (from left to right) form the vertices of the
dominating set.  However, for the graph on the right, we make pebbling
moves so that the second and fourth vertices are selected to be the vertices of the dominating set.  In some cases, moreover, it takes more than the minimum dominating set of vertices to form the minimal domination cover solution.  For example, in Figure 2 below we consider the case of the binary tree with height two, where the minimum dominating set has two vertices, but the minimal dominating set when creating a domination cover solution has three vertices.  This corresponds to several possible starting configurations, e.g., the configuration as pictured; or one with a pebble at the leftmost bottom level vertex and 4 pebbles at the root; or one with 1 and 10 pebbles at the leftmost and rightmost bottom level vertices respectively.
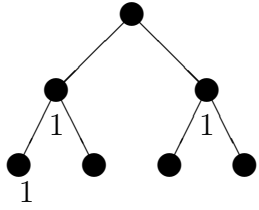
\begin{figure}[tree]
\unitlength 1mm
\begin{center}
\begin{picture}(30,25)
\put(15,25){\circle*{3}}
\put(15,25){\line(-1,-1){11}}
\put(15,25){\line(1,-1){11}}
\multiput(5,15)(20,0){2}{\circle*{3}}
\put(5,15){\line(-1,-2){5}}
\put(5,15){\line(1,-2){5}}
\put(25,15){\line(-1,-2){5}}
\put(25,15){\line(1,-2){5}}
\multiput(0,5)(10,0){2}{\circle*{3}}
\multiput(20,5)(10,0){2}{\circle*{3}}
\put(0,0){1} \put(4,9){1} \put(24,9){1}
\end{picture}
\end{center}
\caption{A reachable minimal configuration of pebbles on $B_2$ that forces
a domination cover solution.} \label{ex1}
\end{figure}

The above two facts constitute the main reason why domination cover
pebbling is nontrivial.  We refer the reader to
\cite{haynes} for additional exposition on domination in graphs.

The paper
is organized as follows.  First, we present some basic results about
domination cover pebbling. The remainder of the paper will consist of
proofs that determine the domination cover pebbling number of $P_n$, the
path on $n$ vertices, $C_n$, the cycle graph on $n$ vertices, and $B_n$,
the complete binary tree with height $n$.  A companion paper by Watson and Yerger (\cite{dcp2})
examines the relation between $\psi(G)$ and structural properties of the graph.

\section{Preliminary Results} We begin by determining the domination cover
pebbling number for some families of graphs.
\begin{thm}
For the complete graph $K_n$ on $n$ vertices, $\psi(K_n) = 1$.
\end{thm}
\noindent This result is obvious since placing a pebble on any vertex dominates
$K_n$.

\begin{thm}
For $ s_1 \geq s_2 \geq \cdots \geq s_r$, let $K_{c_1, c_2, \ldots, c_r}$
be the complete $r$-partite graph with $s_1,s_2,\ldots,s_r$ vertices in
vertex classes $c_1, c_2, \ldots, c_r$ respectively.  Then, for $s_1 \geq
3$, $\psi(K_{c_1, c_2, \ldots, c_r}) = s_1$. If $s_1=2$, $\psi(K_{c_1, c_2, \ldots, c_r}) = 3$.
\end{thm}

\medskip

\noindent{\bf Proof}
First, the configuration with one pebble on all but one of the vertices in $c_1$ does not produce a domination cover solution. So
$\psi(K_{c_1, c_2, \ldots, c_r})$ $> s_1 -1$.  Notice next that if there are
pebbles on vertices in two different $c_i's$, the configuration contains a domination cover
pebbling.  Thus, any pair of pebbles on a vertex along  with another pebbled
vertex can force a domination cover pebbling. So if there are $s_1$
pebbles, the only configuration that has not been considered is the one with one pebble on every vertex in a vertex class that contains $s_i$
 vertices, but this also forces a domination cover pebbling.  Hence, $\psi(K_{c_1, c_2, \ldots, c_r})$ $ = s_1$; notice how we used the condition $s_1\ge3$.
\hfill\qed

For the next theorem we define the wheel graph, denoted $W_n$, to be the
graph with $V(W_n) = {h, v_1,v_2,\ldots,v_n}$, where $h$ is called the hub
of $W_n$, and $E(W_n) = C_n \cup \{hv_1,hv_2,\ldots,hv_n\}$, where $C_n$ denotes the cycle graph on $n$ vertices.
\begin{thm}
For $n \geq 3$, $\psi(W_n) = n-2$.
\end{thm}

\medskip

\noindent{\bf Proof}

First, $\psi(W_n) > n-3$ because placing one pebble on each of $n-3$
consecutive outer vertices leaves a vertex of $W_n$ undominated. If there
is a pair of pebbles on any vertex, move it to the center, and the
domination is complete. Likewise, if there is a pebble in the hub vertex,
$W_n$ is dominated. Thus, consider all configurations containing pebbled
vertices that each contain only one pebble.  If there are $n-2$ vertices
containing pebbles, the two non-pebbled outer vertices are forced to be
dominated since there are only $3$ vertices in all of $W_n$ that contain no
pebbles. Therefore, $\psi(W_n) = n-2$.
\hfill\qed
\section{Domination Cover Pebbling for Paths}
\begin{thm}
\[
\psi(P_n) = \etapn \textrm{,}
\]
for $n \geq 3$, where $n-2 = \alpha_n +
3k_n\equiv \alpha_n\pmod 3$.
\end{thm}

\medskip

\noindent{\bf Proof}
Let $V = V(P_n) = \{v_1,v_2, \ldots ,v_n\}$ with $E(P_n) = \{v_1v_2, \ldots
v_{n-1}v_{n}\}$. Consider the configuration where all pebbles are placed on
$v_1$. We need at least $2^{n-2}$ pebbles to dominate $v_n$. Likewise, we
need at least $2^{n-2} + 2^{n-5} + 2^{n-8} + \cdots + 2^{\alpha_n}$ pebbles to
dominate $\{v_{n},v_{n-1},\cdots , v_{\alpha_{n+1}}\}$. If $\alpha_n = 0$
or $1$, then we have already dominated $P_n$. Otherwise, $\alpha_n = 2$ and
we need one more pebble on either $v_1$ or $v_2$ to dominate $P_n$. Thus,
under this configuration,
 \begin{eqnarray*}
 \psi(P_n) &\geq& 2^{n-2}\sum_{i=0}^{k_n}\frac{1}{8^i} + \floor{{{\frac{\alpha_n}{2}}}} \\
    &=& \etapn,
\end{eqnarray*}
since $\lfloor{\frac{\alpha_n}{2}}\rfloor = 0$ for $\alpha_n = 0$ or $1$ and $\lfloor{\frac{\alpha_n}{2}}\rfloor = 1$
for $\alpha_n=2$.

We now use induction to show that
\[
\psi(P_n) \leq \etapn\textrm{.}
\]
The assertion is clear for $n=3$. Therefore, we assume it is true for all
$P_m$, where $3 \leq m \leq n-1$. Consider an arbitrary configuration of
$P_n$ having $\letapn$ pebbles. Clearly, we can cover dominate
$\{v_{n-2},v_{n-1},v_n\}$ in a finite number of moves with $2^{n-2}$
pebbles or less. Thus, we need to dominate $P_{n-3}$ with the remaining
\[
\etapn - 2^{n-2} = 2^{(n-3)+1} \left(\frac{1-8^{-(k_{n-3}+1)}}{7}\right)
+ \floor{{\frac{\alpha_{n-3}}{2}}}
\]
pebbles, since $\forall n$, $k_n = k_{n-3}+1$ and $\alpha_n=\alpha_{n-3}$.
This number of pebbles is enough to dominate $P_{n-3}$ by hypothesis. Thus,
\[\psi(P_n) \leq \etapn\textrm{,}\]
completing the proof.
\hfill\qed
\section{Domination Cover Pebbling for Cycles} We begin by
proving that placing all the pebbles on one vertex is the ``worst case" configuration that determines the domination
cover pebbling number.
\begin{lem}The value of $\psi(C_n)$ is attained when the original configuration consists of placing
all the pebbles on a single vertex.
\end{lem}

\medskip

\noindent {\bf Proof}\
The proof relies strongly on the fact that the underlying graph is $C_n$, and is by contradiction.  Assume first that the worst configuration consists of more than one set of consecutively pebbled vertices (``islands").  Consider the cardinality of any such set.  This must be at most two, for, were it to be three or more, one could move the pebbles to the inner one or two vertices, thereby causing a larger number of pebbles to be needed to dominate -- a contradiction. Thus each such ``island" must have at most two vertices.  Now consider the effect of moving all the pebbles onto a single such island.  Once again one reaches a contradiction since one would now require more pebbles than $\psi(C_n)$ to cover dominate the graph.  So we must then have the worst configuration consisting of a single island of two pebbled vertices.  Clearly, the worst case is placing $\psi-1$ pebbles on one vertex, say $v_1$, and a single pebble on the other vertex, say $v_2$, since it would now cost more pebbles to reach the $v_2$ side of the cycle.  Had, however, all the pebbles been on $v_1$, we would need at least two more pebbles to dominate the other vertex closest to $v_2$, raising a contradiction.  The statement follows.
\hfill\qed

Since placing all the pebbles on a single vertex is the worst case, we now
determine the value of $\psi(C_n)$.
\begin{thm} Let $C_n$ be a cycle on $n$
vertices. If $n=2m-1,m\ge2$,
\[
    \psi(C_n) = 2^{m+2} \left({\frac{1-8^{-(k_m+1)}}{7}}\right) + \phi_1(m)
\]
and if $n=2m-2,m\ge3$,
\[
    \psi(C_n) = 2^{m+1} \left({\frac{1-8^{-(k_m+1)}}{7}}\right) + 2^{m} \left({\frac{1-8^{-(k_{m-1}+1)}}{7}}\right) + \phi_2(m)\textrm{,}
\]
where $\phi_1(m) = 2\floor{\alpha_m / 2} - \abs{\alpha_m-1}$, $\phi_2(m) =
\floor{\alpha_m / 2} + \floor{\alpha_{m-1} / 2} -
\abs{\alpha_m-1}\abs{\alpha_{m-1}-1}$, $m-2 \equiv \alpha_m\pmod3$, and
$m-2 = \alpha_m + 3k_m$.
\end{thm}

\medskip

\noindent{\bf Proof}
  By Lemma 5, we assume all $\psi(C_n)$ pebbles are on $v_1 \in
C_n$. If $n=2m-1$, there are two identical $m$ paths to cover. We can cover
these with $2\psi(P_m)$ pebbles. We notice that
$v_1$ may be in both dominating sets; $\abs{\alpha_m-1}=1$ if $v_1$ is double counted. If
$n=2m-2$, there are two paths $P_1, P_2 \in C_n$ with $m-1=\abs{P_2} = \abs{P_1}
-1$. Thus we can cover these two paths with $\psi(P_m) + \psi(P_{m-1})$ pebbles.
Likewise in this case, we may have double-counted vertex $v_1$;
$\abs{\alpha_m-1}\abs{\alpha_{m-1}-1} = 1$ in these cases, i.e. $\alpha_m\equiv0\pmod3; \alpha_{m-1}\equiv2\pmod3$.
Thus we compute the domination cover pebbling number as follows. When
$n=2m-1$
\begin{eqnarray*}
    \psi(C_n) &=& 2\psi(P_m)-\abs{\alpha_m-1} \\
            &=& 2^{m+2} ((1-8^{-(k_m+1)})/7) + 2\floor{\alpha_m / 2} - \abs{\alpha_m-1}
\end{eqnarray*}
and if $n=2m-2$,
\begin{eqnarray*}
    \psi(C_n) &=& \psi(P_m)+ \psi(P_{m-1}) - \abs{\alpha_m-1}\abs{\alpha_{m-1}-1} \\
            &=& 2^{m+1} ((1-8^{-(k_m+1)})/7) + 2^{m}((1-8^{-(k_{m-1}+1)})/7) \\
            && + \floor{\alpha_m / 2} + \floor{\alpha_{m-1} / 2} -
            \abs{\alpha_m-1}\abs{\alpha_{m-1}-1},
\end{eqnarray*}
as asserted.
\hfill\qed

\section{Binary Trees}

In this section, we will compute the domination cover pebbling number for
the family of complete binary trees.
Recall that a complete binary tree, denoted by $B_n$, is a tree of height
$n$, with $2^i$ vertices at distance $i$ from the root.  Each vertex of
$B_n$ has two ``children", except for the set of $2^n$ vertices that are
distance $n$ away from the root, none of which have children.  The root will be denoted by $\rho=\rho_n$

\begin{thm} $\psi(B_0)=1, \psi(B_1)=2, \psi(B_2)=11, \psi(B_3)=81, \psi(B_4)=609.$
\end{thm}

\medskip

\noindent{\bf Proof}
The fact that $\psi(B_0)=1$ and $ \psi(B_1)=2$ are obvious.  We next show that $\psi(B_2)=11$, as predicted by the general formula of Theorem 8.  In Figure 3, we exhibit a configuration
of $10$ pebbles on $B_2$ that does not force a domination cover solution.

\begin{figure}[htb]
\unitlength 1mm
\begin{center}
\begin{picture}(30,25)
\put(15,25){\circle*{3}}
\put(15,25){\line(-1,-1){11}} 
\put(15,25){\line(1,-1){11}}
\multiput(5,15)(20,0){2}{\circle*{3}} 
\put(5,15){\line(-1,-2){5}} 
\put(5,15){\line(1,-2){5}}
\put(25,15){\line(-1,-2){5}} 
\put(25,15){\line(1,-2){5}}
\multiput(0,5)(10,0){2}{\circle*{3}} 
\multiput(20,5)(10,0){2}{\circle*{3}} 
\put(0,0){1}  \put(30,0){9}
\end{picture}
\end{center}
\caption{A configuration of 10 pebbles on $B_2$ that does not force a
domination cover solution.} \label{ex2}
\end{figure}
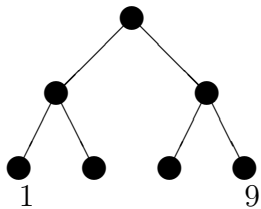

We will now show that $\psi(B_2) \leq 11$.  Arbitrarily
place $11$ pebbles on $B_2$.  Consider the following three subcases based on
the number of pebbles on each of the two $B_1$'s connected to the root of
$B_2$.

Case $1$:  Suppose there are at least two pebbles on each of the two
$B_1$'s.  It takes at most two pebbles for each of the $B_1$'s to be
dominated.  Hence, after dominating each of the $B_1$'s there are seven
pebbles left.  If there is a pebble on either of the two root vertices of
the two disjoint copies of $B_1$, then we have dominated the root of $B_2$.
Otherwise, it is always possible to move a pebble to the root of one of the $B_1$'s, thus dominating the root
$\rho_2$. This process induces a domination cover solution
of $B_2$, completing this case.

Case $2$:  Suppose that neither $B_1$ contains two or more pebbles.  Then
there are at least $9$ pebbles on the root of $B_2$.  Pebble the root of
each of the $B_1$'s, and this case is complete.

Case $3$:  Suppose that one copy of $B_1$ contains two or more pebbles,
call it $B_1^*$, and the other copy does not.  Then all of the pebbles on
$B_1^*$ except for two can be used, together with any pebbles already on $\rho_2$, to place two pebbles on $\rho_2$, which can be used to dominate the other $B_1$.

We now show that $\psi(B_3) = 81$.  First, we have constructed
in Figure 4 a configuration of $80$ pebbles that does not produce a domination
cover pebbling.

\begin{figure}[htb]
\begin{center}
\unitlength 1mm
\begin{picture}(80,40)
\put(37.5,35){\circle*{3}}
\put(37.5,35){\line(-5,-2){22}} 
\put(37.5,35){\line(5,-2){22}}
\multiput(15,25)(45,0){2}{\circle*{3}} 
\put(15,25){\line(-1,-1){11}} 
\put(15,25){\line(1,-1){11}}
\put(60,25){\line(-1,-1){11}} 
\put(60,25){\line(1,-1){11}}
\multiput(5,15)(20,0){2}{\circle*{3}} 
\multiput(50,15)(20,0){2}{\circle*{3}} 
\put(5,15){\line(-1,-2){5}} 
\put(5,15){\line(1,-2){5}}
\put(25,15){\line(-1,-2){5}} 
\put(25,15){\line(1,-2){5}}
\put(50,15){\line(-1,-2){5}} 
\put(50,15){\line(1,-2){5}}
\put(70,15){\line(-1,-2){5}} 
\put(70,15){\line(1,-2){5}}
\multiput(0,5)(10,0){2}{\circle*{3}} 
\multiput(20,5)(10,0){2}{\circle*{3}} 
\multiput(45,5)(10,0){2}{\circle*{3}} 
\multiput(65,5)(10,0){2}{\circle*{3}} 
\put(0,0){1} \put(20,0){1} \put(45,0){1} \put(75,0){77}

\end{picture}
\caption{A configuration of 80 pebbles on $B_3$ that does not force a
domination cover solution.} \label{ex3}
\end{center}
\end{figure}
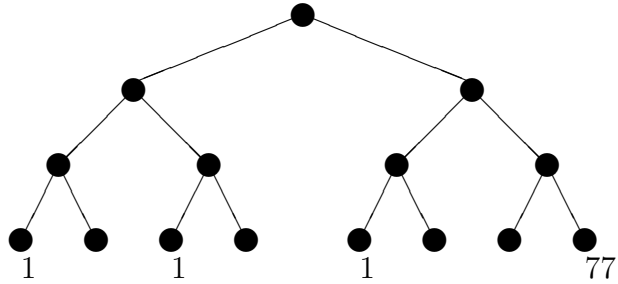

Now suppose that we are given a configuration of $81$ pebbles on $B_3$.  We
wish to force a domination cover solution on $B_3$.  If there are fewer
than $11$ pebbles on each of the two disjoint $B_2$ subtrees in $B_3$, then
we can use 17 of the $61+$ pebbles on the root vertex to produce a domination
cover solution.  If there are at least $11$ pebbles on both of the disjoint
$B_2$ subtrees, then we can dominate the
root vertex with the $59$ remaining pebbles as follows:  One subtree must contain 30 of these ``extra" pebbles, and thus by the pigeonhole principle one of the 4 paths leading from $\rho_3$, the root of $B_3$, to the bottom of the tree must have 8 pebbles on it, enough to send a pebble to the root, since the pebbling number $\pi(P_4)$ of the 4-path (3 edges) is 8.

Next, consider the case when only one of the two disjoint $B_2$ subgraphs,
call it $B_2^*$, contains at least $11$ pebbles.  There are at most $70$
pebbles somewhere on the graph that must now be succesfully used to dominate the other
$B_2$, call it $B_2'$, and $\rho_3$.  Our strategy will be to move as many pebbles as possible from $B_2^*$ to $\rho_3$ while still leaving $B_2^*$ dominated.  The pebbles placed on $\rho_3$ in this fashion will then be used to reach the two ``middle" vertices in $B_2'$ so as to dominate it.  Notice that two single pebbles
on the bottom row of $B_2'$, each with a different parent, does not
decrease the number of pebbles required to dominate $B_2'$ when there are no other pebbles on $B_2'$ and all pebbles used in the domination emanate at the root. Also, any preexisting pebbles at $\rho_3$ [or any pebbles on $B_2'$ other than the above-mentioned two] only make our strategy easier to implement, so assume that 68 of the extraneous pebbles are on $B_2^*$ and the other two on the bottom row of $B_2'$ as specified above.

Call the vertices of $B_2^*$ $a$ (its root); $b_1,b_2$ (the ``middle vertices"); and $c_1,c_2,c_3,c_4$ (the bottom vertices).  We thus need to force 9 pebbles onto $\rho_3$, failing which we need 8 pebbles on $\rho_3$ and a pebble on $a$.  Variations of the argument that follows will be used throughout this paper.  In order to accomplish our task, we will use (in addition to the 68 extraneous pebbles)  the 11 pebbles ``reserved" to dominate $B_2^*$. Now each pebble sent to $\rho_3$ causes a net reduction of at most 8 pebbles and a pebble sent to $a$ causes at most 4 pebbles to be lost.  Since $9\times8=72>68$, it appears that we can't always send 9 pebbles to $\rho_3$, so let's try to send 8 to the root and one to $a$.  We claim and prove next that {\it we can send a pebble to $\rho_3$ as long as $B_2^*$ has a total of 18 or more pebbles on it}.  If there are 18 pebbles on $B_2^*$, one $``a-b-c"$ path must contain 5 pebbles, say $a-b_1-c_1$.  The maximum possible number of pebbles on this path is 7, or else we could send a pebble to the $\rho_3$. If there are exactly 5 pebbles on $a-b_1-c_1$ (13 left over), one of the remaining 4 vertices must have 4 pebbles on it, one of which can reach $a$.  Another pebble can be put on $a$ using the 5 pebbles on the $a-b-c$ path.  We can reach $\rho_3$.  If there are 6 pebbles on the $a-b-c$ path, there is one vertex with 3 pebbles.  This must be $c_3$ or $c_4$, say $c_3$, which must, moreover, have {\it exactly} 3 pebbles (otherwise 2 pebbles can be placed on $a$.)  The next vertex guaranteed to have 3 pebbles on it by the pigeonhole principle can be checked to lead to $\rho_3$ being pebbleable.  Finally, suppose the $a-b-c$ path has  7 pebbles on it.  The vertex guaranteed to have (exactly) 3 pebbles on it must be (WLOG) $c_3$.  But there is another such vertex, which can again be checked to lead to $\rho_3$ being reached.

We had started with 68+11 pebbles on $B_2^*$.  Thus 8 pebbles can be sent to $\rho_3$ with as few as 15 left on $B_2^*$.  But 15 pebbles do imply that there is an $a-b-c$ path with 4 pebbles, enough to reach $a$. Thus, $\psi(B_3) =
81$.

We are ready to prove that $\psi(B_4)=609$.  Starting on the left, place one pebble at every alternate vertex on the bottom row with one exception:  For the last two vertices, we place no pebbles on the left vertex and 601 pebbles on the right vertex.  This construction, similar to those used for $n=2,3$ will be a canonical one that we will use in general later. The most efficient pebble domination would be to place a pebble at $\rho_4$ and one on each vertex at the next to bottom level of the tree.  It is elementary to check that this cannot be done.  Hence $\psi(B_4)\ge609$.  We now prove that $\psi(B_4)\le609$.

If there are fewer than 81 pebbles on each 3-subtree, $\rho_4$ must have have at least 449 pebbles, but only 64 of these are required to pebble the 8 vertices in the next to bottom row.  This, together with one more pebble at $\rho_4$, completes the pebble domination.  If there are 81+ pebbles on each 3-subtree, we use 81 to dominate each subtree, leaving us with 447 pebbles to dominate $\rho_4$.  224 of these must be on one subtree, so one of the 8 paths leading from $\rho_4$ to the bottom of this subtree must have 28 pebbles, enough to reach $\rho_4$.  As before the most complicated case is when one subtree, $B_3^*$ has 81+ pebbles and the other doesn't.  We employ the same pebbling strategy as for $n=3$.  Put one pebble on each of 4 bottom vertices in $B_3'$ (no two of whom share a parent).  Assume that $\rho_4$ is unpebbled, and that there are just the four aforementioned pebbles on $B_4'$.  We thus have 605 pebbles on $B_3^*$.  Our proof deviates here from the $n=3$ case, and in general, the case of $n\equiv1\pmod3$ will be seen later to be trickier than the rest.  We can clearly assume, given the strategy being used, that the root $\rho_3^*$ of $B_3^*$ is unpebbled.   It is a straightforward calculation (similar to the $n=3$ case) to verify that $\psi(B_3^*\setminus\rho_3^*)$, the domination pebbling number of $B_3^*$ minus its root, equals 77.  We will now seek to place 33 pebbles on $\rho_4$ -- adequate to pebble the next to bottom row of $B_3'$ while leaving a pebble at $\rho_4$ -- while never dropping below 77 pebbles on $(B_3^*\setminus\rho_3^*)$. The extra pebble on $\rho_4$ will serve to dominate $\rho_3^*$.  We have 605 pebbles on $(B_3^*\setminus\rho_3^*)$.  We claim that as long as the number of pebbles does not drop to below 77, we can get a pebble to $\rho_4$ at a loss of $\le16$ pebbles.  Since $(33\times16)+77=605$, we will be done if we can show that 93 (=77+16) pebbles on $(B_3^*\setminus\rho_3^*)$ suffice to send a pebble to $\rho_4$.  The pigeonhole principle, the fact that $\rho_4$ has no pebbles on it, and the observation that no path from $\rho_4$ to the bottom can have 16+ pebbles, implies that we must have two $a-b-c$ paths with 12+ pebbles each.  If these paths are disjoint, we can move two pebbles to $\rho_3^*$ and thus one to $\rho_4$, so the paths must overlap.  There are two possibilities:  the paths may be of the form $a-b_1-c_1$ and $a-b_1-c_2$, or they may overlap in just the vertex $a$.  In the first case, the worst case scenario is, e.g., when there are no pebbles on $a$, one on $b_1$, and 11 on each of the $c$'s.  The second case is similar.  But these configurations too force a pebble onto $\rho_4$.  Thus $\psi(B_4)=609$.
\hfill\qed
\vfill\eject
\begin{thm} For $n \geq 2$,
\begin{eqnarray*} &&{}\psi(B_n)\\ &&{}= (2^{n-1}-1) +
\sum_{i=0}^{\lfloor\frac{n-1}{3} \rfloor} \lr2^{3i + 1} + \sum_{j=1}^{n-3i-2}
2^{j-1} 2^{3i+2j+1}\rr  + \sum_{k=1}^{\lfloor \frac{n+1}{3} \rfloor} 2^{n-3k
+ 1} 2^{2n-3k+2} + \gamma_n,\\
&&{}=T_{1,n}+T_{2,n}+T_{3,n}+T_{4,n}\quad{\rm say,}\end{eqnarray*}
where $T_{i,n}$ denotes the $i^{\rm th}$ term in the above sum, and $\gamma=\gamma_n=2^{n-1}$ if $n \equiv 0 \pmod 3$, and $\gamma = 0$ otherwise.
\end{thm}

\medskip

\noindent{\bf Proof}  First we will prove that for $n\ge 2$\begin{eqnarray*}
\psi(B_n) & > & (2^{n-1}-1) +
\sum_{i=0}^{\lfloor\frac{n-1}{3} \rfloor} [2^{3i + 1} + \sum_{j=1}^{n-3i-2}
2^{j-1} 2^{3i+2j+1}] \\
& & + \sum_{k=1}^{\lfloor \frac{n+1}{3} \rfloor}
2^{n-3k + 1} 2^{2n-3k+2} + \gamma-1.
\end{eqnarray*}
Consider the following initial configuration of pebbles that generalizes that in Figures 3 and 4.  Starting at the left, place one pebble at
each of $2^{n-1} -1$ vertices on the bottom row, no two of which share a parent.  This leaves the rightmost two vertices on the bottom row unpebbled; we place all the remaining pebbles on one of these vertices, denoted by $v$.
We will endeavor to pebble entire rows in the most efficient way -- the rows to be pebbled are specified by working upwards from the bottom of both subtrees -- and we make pebbling moves from $v$, if possible, so that one
pebble is placed on every vertex in every third row, starting, in each subtree, with the row
that is next to the bottom row. This is clearly the best strategy, since it ``costs" the most to reach the bottom row of the left subtree.  Note that the $\gamma$ term enters iff $n\equiv0\pmod3$, since then the three topmost vertices of $B_n$ would be left unpebbled, but we can complete the domination pebbling of the graph by pebbling the root of the right subtree.  This requires $\gamma$ pebbles.    If we consider rows as single vertices,
this would be analogous to the configuration of pebbles required to get an
optimal domination cover pebbling bound for $P_n$, except that we pebble every third vertex  counted from {\it both} ends, adding a ``central $\gamma$-correction" if needed.

We start by stating that in order to find a domination cover solution
for the subtree that is on the other side of the root vertex $v$ it takes
$$\sum_{k=1}^{\lfloor \frac{n+1}{3} \rfloor} 2^{n-3k + 1} 2^{2n-3k+2}$$
pebbles as follows:  Consider the next to bottom row.  There are
$2^{n-2}$ vertices that must have a pebble placed on them.  For each
vertex, it takes $2^{2n -1}$ pebbles from vertex $v$, for a total of
$2^{3n-3}$ pebbles.  This is the number of pebbles counted in the $k=1$
term of the sum, since $2^{n-3+1}2^{2n - 3 + 2}= 2^{3n-3}$.  We leave the rest of the details -- of verifying that the stated expression $\sum_{k=1}^{\lfloor \frac{n+1}{3} \rfloor} 2^{n-3k + 1} 2^{2n-3k+2}$ does indeed represent the above pebbling process of the left subtree -- to the reader.   A similar (but somewhat more complicated) computation can be performed to verify the first sum represents the pebbling of the subtree on the same side as $v$ in the manner desired -- {\it except that, due to the last (-1) term, the vertex right above the ``pebble source" remains unpebbled.}  This configuration leaves the other sibling of $v$ undominated.  This proves the claim.

We now proceed to prove the bound by induction.  First we note that a simple Maple computation reveals that the values given by the putative formula for $\psi(B_n)$ are, for $n=2,\ldots10$, equal to 11, 81, 609, 4777, 38105, 304473, 2434969, 19478809, 155827481; the first three of these have already been proved to be correct in Theorem 7. Note that the asymptotic ratio of the terms appears to be converging to 8 rapidly, reflecting the fact that the dominant term in $\psi(B_n)$ is $2^{3n-3}$.
Suppose that the value of
$\psi(B_{n-1})$ is as stated in the theorem for $n\ge5$, i.e. $n-1\ge4$ (our induction will only work for these cases). Place $\psi(B_{n})$ pebbles on $B_n$.  As before, we will consider three cases depending
upon whether there are enough pebbles in each of the two disjoint copies of
$B_{n-1}$ connected to the root of $B_n$.

First, suppose that neither copy
contains $\psi(B_{n-1})$ pebbles. In this case, there are at least $\psi(B_n)-2\psi(B_{n-1})+1$ pebbles on the root.  We claim that this number is at least as large as $4\psi(B_{n-1})+1$, a number that would allow us to move $\psi(B_{n-1})$ pebbles onto the root of each subtree while retaining one pebble on the root, thus completing the domination of $B_n$. It suffices to show that
\begin{equation}\psi(B_n)\ge6\psi(B_{n-1})\end{equation} in order for the above to be true.  We have
\begin{equation}\psi(B_n)\ge2^{3n-3}\end{equation}
 by considering only the $k=1$ term of $T_{3,n}$.  Also,
\begin{equation}6(T_{1,n-1}+T_{4,n-1})\le 3\cdot2^{n};\end{equation}
\begin{equation}
6T_{3,n-1}=6\sum_{k=1}^{\floor{\frac{n}{3}}}2^{3(n-1)-6k+3}\le\frac{6\cdot2^{3n}}{63};
\end{equation}
and
\begin{eqnarray}
6T_{2,n-1}&=&6\sum_{i=0}^{{\floor{\frac{n-2}{3}}}}2^{3i+1}+6\sum_{i=0}^{{\floor{\frac{n-2}{3}}}}2^{3i}\sum_{j=1}^{n-3i-3}2^{3j}\nonumber\\
&\le&\frac{24}{7}2^n+6\sum_{i=0}^{{\floor{\frac{n-2}{3}}}}2^{3i}\cdot\frac{8(8^{n-3i-3})}{7}\nonumber\\
&\le&\frac{24}{7}2^n+\frac{3\cdot8^n\sum_{i=0}^{{\floor{\frac{n-2}{3}}}}8^{-2i}}{224}\nonumber\\
&\le&\frac{24}{7}2^n+\frac{6}{441}2^{3n}.
\end{eqnarray}
Equations (2) through (5) show that (1) holds if
$$2^{3n-3}\ge3\cdot2^n+\frac{6\cdot2^{3n}}{63}+\frac{24}{7}2^n+\frac{6}{441}2^{3n}$$
or if
$$\frac{1}{8}-\frac{6}{63}-\frac{6}{441}\ge\frac{45}{7}2^{-2n},$$
which holds for $n\ge5$.  Of course the fact that $\psi(B_n)\ge6\psi(B_{n-1})$ holds for $n=3,4$ as well.

Next, suppose that both copies contain at least $\psi(B_{n-1})$ pebbles. In
this case we can use $\psi(B_{n-1})$ pebbles to construct a domination
cover pebbling for each subtree.  At least one subtree thus has $2\psi(B_{n-1})\ge\frac{2^{3n+1}}{64}$ extra pebbles (recalling that $\psi(B_n)\ge6\psi(B_{n-1})$ and $\psi(B_{n-1})\ge2^{3n-6}$), so at least one of the $2^{n-1}$ $n+1$-paths leading to the root vertex (from the bottom of the subtree) has at least $\frac{2^{2n}}{16}$ pebbles.  For $n\ge4$ this number exceeds $2^n$, the (regular) pebbling number of $P_{n+1}$.  We can thus reach the root vertex.  [If $n=3$, the exact values of $\psi(B_3), \psi(B_2)$ show that we must have 30 pebbles on one of the 2-subtrees and we can reach the root as desired since one of the paths to the root must have eight pebbles.]

Finally, suppose that only one copy of $B_{n-1}$, call it $B_{n-1}^*$, contains
at least $\psi(B_{n-1})$ pebbles.  We need to do a more careful analysis, since the strategy to be employed (as in the small cases studied in Theorem 7) is to move all extraneous pebbles in $B_{n-1}^*$ to the root of the tree and then cover pebble dominate the other subtree, say $B_{n-1}'$, from the root using an ``every third row" dominating set.  Note that any  pebbles in
$B_{n-1}'$ can substitute for at least one pebble on the root vertex.  We may thus  assume the worst case scenario in which {\it all} the $\psi(B_{n})$ pebbles are in $B_{n-1}^*$ -- except for the $2^{n-2}$ non-siblings in the bottom  row which can each be assumed to have a single pebble on them, since this does not lessen the pebbling number for the left subtree starting at the root.

We start by showing that we may place a pebble on $\rho_n$ (at a loss of at most $2^{n}$ pebbles) if there are at least $2^{3n-6}$ pebbles on $B_{n-1}^*$ and thus, since
$$\psi(B_{n-1}^*)\ge 2^{3n-6}+2^{n-2}-1+\sum_{i=0}^{\floor{\frac{n-2}{3}}}2^{3i+1}\ge2^{3n-6}+2^{n-2},$$ if there are at least
$\psi(B_{n-1}^*)-2^{n-2}$ pebbles on $B_{n-1}^*$  (Of course we will cause problems if we send a pebble to the root when there are fewer than $\psi(B_{n-1}^*)+2^n$ pebbles on $B_{n-1}^*$.)  To see this, note that one of the $2^{n-1}$ paths leading to the root from the bottom of $B_{n-1}^*$ must have $2^{3n-6}/2^{n-1}=2^{2n-5}$ pebbles on it.  This number exceeds $2^n$ if $n\ge5$, so we can send a pebble to the root.  Similarly, a pebble may be sent to the root of $B_{n-1}^*$ if there are at least $2^{3n-6}$ pebbles on $B_{n-1}^*$.  This costs at most $2^{n-1}$ pebbles.
We have $\psi(B_n)-2^{n-2}$ pebbles on $B_{n-1}^*$, and $\psi(B_n)-\psi(B_{n-1})-2^{n-2}$ of these are available  to pebble $B_{n-1}'$.  Let's compute $D(n) = \psi(B_n) -
\psi(B_{n-1})-2^{n-2}$; by the above argument, we will, e.g. be able to send $\floor{D(n)/2^n}$ pebbles to the root.  We get:
\begin{eqnarray*}
D(n)& = &
 2^{n-1} + \sum_{i=0}^{\lfloor\frac{n-1}{3} \rfloor} [2^{3i + 1} +
\sum_{j=1}^{n-3i-2} 2^{3i+3j}] + \sum_{k=1}^{\lfloor
\frac{n+1}{3} \rfloor} 2^{3n-6k+3} + \gamma_n\\ & &   - \left[
2^{n-2} + \sum_{i=0}^{\lfloor\frac{n-2}{3} \rfloor} [2^{3i + 1} +
\sum_{j=1}^{n-3i-3} 2^{3i+3j}]  + \sum_{k=1}^{\lfloor \frac{n}{3}
\rfloor} 2^{3n-6k} +\gamma_{n-1}\right]-2^{n-2}
\end{eqnarray*}
\begin{eqnarray*}
& \geq & \sum_{i=0}^{\lfloor \frac{n-1}{3}\rfloor}2^{3i+1}-\sum_{i=0}^{\lfloor \frac{n-2}{3}\rfloor}2^{3i+1}+\sum_{i=0}^{\lfloor \frac{n-1}{3}\rfloor} 2^{3n-6i-6} + \frac{7}{8}
\sum_{k=1}^{\lfloor \frac{n+1}{3} \rfloor} 2^{3n-6k+3} + \gamma_n-\gamma_{n-1} \\
& = & \sum_{i=0}^{\lfloor \frac{n-1}{3}\rfloor}2^{3i+1}-\sum_{i=0}^{\lfloor \frac{n-2}{3}\rfloor}2^{3i+1}+\sum_{i=1}^{\lfloor \frac{n-1}{3}\rfloor + 1} 2^{3n-6i} +
\frac{7}{8} \sum_{k=1}^{\lfloor \frac{n+1}{3} \rfloor} 2^{3n-6k + 3} + \gamma_n-\gamma_{n-1}\\
& \geq & \sum_{i=0}^{\lfloor \frac{n-1}{3}\rfloor}2^{3i+1}-\sum_{i=0}^{\lfloor \frac{n-2}{3}\rfloor}2^{3i+1}+\sum_{k=1}^{\lfloor \frac{n+1}{3} \rfloor} 2^{3n-6k +
3} + \gamma_n-\gamma_{n-1}.
\end{eqnarray*}
{\it Case 1}.  If $n\equiv0\pmod3$, we cover dominate $B_{n-1}^*$ using $\psi(B_{n-1})$ pebbles.  Note that $\gamma_n-\gamma_{n-1}=2^{n-1}$ and these $2^{n-1}$ pebbles are used to pebble the root of $B_{n-1}^*$, so that $\rho_n$ is dominated in addition to $B_{n-1}^*$.  Also, $\sum_{i=0}^{\lfloor \frac{n-1}{3}\rfloor}2^{3i+1}-\sum_{i=0}^{\lfloor \frac{n-2}{3}\rfloor}2^{3i+1}=0$.  It follows that $\sum_{k=1}^{\lfloor \frac{n+1}{3} \rfloor} 2^{2n-6k +
3}$ pebbles can be placed on the root, and it is easy to check that these suffice to pebble dominate $B_{n-1}'$ by placing pebbles on every third row, starting with the next to bottom row.

\noindent{\it Case 2}.  If $n\equiv2\pmod3$, $\gamma_n-\gamma_{n-1}$ and $\sum_{i=0}^{\lfloor \frac{n-1}{3}\rfloor}2^{3i+1}-\sum_{i=0}^{\lfloor \frac{n-2}{3}\rfloor}2^{3i+1}$ both equal zero.  We send $\sum_{k=1}^{\lfloor \frac{n+1}{3} \rfloor} 2^{2n-6k +
3}$ pebbles to $\rho_n$.  These suffice to dominate the left subtree, and in particular the root of
$B_{n-1}'$ is pebbled, so that $\rho_n$ is dominated.

\noindent{\it Case 3}. $n\equiv1\pmod3$.  This is the delicate case.  First, $\sum_{k=1}^{\lfloor \frac{n+1}{3} \rfloor} 2^{2n-6k +
3}$ pebbles are sent to $\rho_n$.  However,
\begin{equation}\gamma_n-\gamma_{n-1}+\sum_{i=0}^{\lfloor \frac{n-1}{3}\rfloor}2^{3i+1}-\sum_{i=0}^{\lfloor \frac{n-2}{3}\rfloor}2^{3i+1}=-2^{n-2}+2^{n}=\frac{3}{4}2^n,\end{equation}
insufficient to place a crucial extra pebble on $\rho_n$.  This is needed since the $\sum_{k=1}^{\lfloor \frac{n+1}{3} \rfloor} 2^{2n-6k +
3}$ pebbles on the root are only sufficient to pebble dominate $(B_{n-1}'\setminus\rho_{n-1}')$.  The root $\rho_{n-1}'$ of the left subtree could have been dominated by an extra pebble on $\rho_n$.  However such a pebble would cause the root of $B_{n-1}^*$ to be possibly double dominated, which seems sub-optimal.  We resolve the problem by mimicking the $n=4$ case. The induction from $n-1$ to $n\equiv1\pmod3$ proceeds as follows:  We use $\psi(B_{n-1}^*\setminus\rho_{n-1}^*)$ pebbles to dominate $(B_{n-1}^*\setminus\rho_{n-1}^*)$.  It is easy to check, after all the work done above, that $\psi(B_{n-1}^*\setminus\rho_{n-1}^*)=\psi(B_{n-1}^*)-2^{n-2}$.  We have gained the extra $\frac{1}{4}2^n$ pebbles we need, since we do not pebble $\rho_{n-1}^*$.  The modified value in (6) is $2^n$.  $\rho_n$ is thus pebbled, and the roots of both subtrees are dominated as a result.   This completes the proof.
\hfill\qed
\section{Open Problems}  We are confident that this paper, together with the companion paper \cite{dcp2}, will spark interest in the question of domination cover pebbling.  Determination of the $\psi$ values for several other families of graphs is an obvious open question. Teresa Haynes has raised the question of determining the domination cover pebbling number when the pebbles {\it must} reach a {\it miminum} dominating set.  Other open problems are raised in \cite{dcp2}.
\section{Acknowledgments}  This research leading to this paper was conducted during the Summer of 2004 under the supervision of the second-named author -- while the rest of the authors were either completing their Master's thesis in Operator Theory (Teguia) or participating in the ETSU REU (Gardner, Vuong, Watson, Yerger).  All but Teguia were supported by NSF Grant DMS-0139286.  Graduate school affiliations are listed for all but Watson and Godbole, who, respectively, are yet to enter graduate school and last attended graduate school 20 years ago.

\end{document}